\documentclass[12pt
]{amsart}
\usepackage{amssymb}
\usepackage[dvips]{graphics}
\font\bb=msbm10 at9.98pt
\def\semidirect{\hbox{$\;$\bb\char'156$\;$}}

\newcommand{\Z}{{\mathbb Z}}

\newcommand{\C}{{\mathbb C}}
\newcommand{\Q}{{\mathbb Q}}

\newcommand{\dontprint}[1]
{\relax}
\newtheorem%
{thm}{Theorem}[section]
\newtheorem%
{proposition}[thm]{Proposition}
\newtheorem%
{lemma}[thm]{Lemma}
\newtheorem%
{lemmadef}[thm]{Lemma-Definition}
\newtheorem%
{corollary}[thm]{Corollary}
\newtheorem%
{conjecture}[thm]{Conjecture}

\title[{}]{Even powers of divisors and elliptic zeta values}
\author[{}]
{Giovanni Felder${}^{1}$ 
\and Alexander Varchenko${}^{2}$}
\thanks{${}^1$Supported in part by the Swiss National Science Foundation and, while at MSRI, by
NSF grant DMS-9810361} 
\thanks{${}^2$Supported in part by NSF grant  DMS-9801582}
\address{G. F.: MSRI, 1000 Centennial Drive, Berkeley, CA 94720, USA
and
Departement of Mathematics,
  ETH-Zentrum, 8092 Z\"urich, Switzerland}
\email{felder@math.ethz.ch}
\address{A. V.: Department of Mathematics, University of
  North Carolina at Chapel Hill, Chapel Hill, NC
  27599-3250, USA} 
\email{anv@email.unc.edu}
\date{May 2002}
\begin{document}
\begin{abstract}
We introduce and study {\em elliptic zeta values},
a two-parameter deformation of the values of 
Riemann's zeta function at positive integers. They are essentially
Taylor coefficients of the logarithm of the elliptic gamma function,
and inherit the functional equations of this function. Elliptic
zeta values at even integers are related to Eisenstein series and
thus to sums of odd powers of divisors. The elliptic zeta values
at odd integers can be expressed in terms of generating series of
sums of even powers of divisors.
\end{abstract}
\maketitle

\newcommand{\Matrix}[4]{\left(\begin{array}{cc}#1&#2\\ #3&#4\end{array}\right)}
\centerline{\it Dedicated to Igor Frenkel on the occasion of his
$50^{\text{th}}$ birthday}

\section{Introduction}
Let $k$ be a positive integer.
The generating function of the sum of $k\!-\!1$st powers of divisors of
positive integers
\[
\sum_{n=1}^\infty\sigma_{k-1}(n)q^n,\qquad
\sigma_{k-1}(n)=\sum_{d\,|\,n}d^{k-1},
\]
converges in $|q|<1$. If $k-1$ is odd,
it has interesting transformation properties 
under the modular
group $\mathrm{SL}(2,\Z)$. 
Indeed, it is closely related to the Eisenstein series
\[
G_k(\tau)=\frac12\sum_{(m,n)\neq(0,0)}(m\tau+n)^{-k}.\]
Let
\[
D_k(q)=\frac{(-2\pi i)^{k}}{(k-1)!}\sum_{n=1}^\infty\sigma_{k-1}(n)q^n.
\]
Then
\[
G_k(\tau)=\zeta(k)+D_{k}(q),
\qquad q=e^{2\pi i\tau},\qquad k=2,4,6,\dots,
\]
where $\zeta(s)=\sum_{n=1}^\infty n^{-s}$ is the
Riemann zeta function. 
For $k\geq 4$, 
these functions are modular forms of weight $k$:
\[
G_k\left(\frac{a\tau+b}{c\tau+d}\right)=(c\tau+d)^{k}G_k(\tau),
\qquad\Matrix abcd\in \mathrm{SL}(2,\Z).
\]
See \cite{Serre} for
proofs of these facts, but notice that there
$G_k$ denotes in our notation $G_{2k}/2$.
As $G_k(-1/\tau)=\tau^kG_k(\tau)$, it follows that, 
if $k$ is even,
\[
\lim_{\tau\to 0}\tau^{k}D_k\left(e^{2\pi i\tau}\right)
=\zeta(k).
\]
In other words, we may regard $\tau^{k}D_k(e^{2\pi i\tau})$
as a one-parameter deformation of the zeta value $\zeta(k)$, and
this is also true if $k$ is odd.
For odd $k$, however, $D_k(q)$ does not have obvious modular
properties. The purpose of this note is to show that $D_k(q)$
can be embedded into a two-parameter deformation of $\zeta(k)$,
the {\em elliptic zeta value} at $k$,
which obeys identities of modular type.
These identities are essentially equivalent to modularity in the even
case but are of a different nature in the odd case. They are related
to (and a consequence of) the three term relations
of the elliptic gamma function \cite{FV}.
\subsection*{Acknowledgments} The first author wishes to thank Daan Krammer for
a useful discussion. We are grateful to Don Zagier for several interesting 
comments and suggestions, and for providing the material of Section \ref{s-5}.

\section{Differences of modular forms}
Let $H$ be the upper half-plane 
$\mathrm{Im}\,\tau>0$.

\begin{proposition} Let $k$ be an even positive integer.
Suppose $Z(\tau,\sigma)$ is a holomorphic
function on $H\times H$ admitting an expansion
\[
Z(\tau,\sigma)=\sum_{n,m=0}^\infty
a_{n,m}q^nr^m,\qquad
q=e^{2\pi i\tau},\quad r=e^{2\pi i\sigma},
\]
and obeying $Z(\tau,\sigma)=-Z(\sigma,\tau)$.
Then the following statements are equivalent:
\begin{enumerate} 
\item[(i)] $Z(\tau,\sigma)=G(\tau)-G(\sigma)$ for
some modular form $G$ of weight $k$.
\item[(ii)] $Z$ obeys the three-term relations
\begin{eqnarray*}
Z(\tau,\sigma)&=&Z(\tau,\tau+\sigma)+Z(\tau+\sigma,\sigma),\\
Z(\tau,\sigma)&=&
\tau^{-k}Z\left(-\frac1\tau,\frac\sigma\tau\right)
+\sigma^{-k}Z\left(-\frac\tau\sigma,-\frac1\sigma\right),
\end{eqnarray*}
for all $\sigma,\tau\in H$ such that $\sigma/\tau\in H$.
\end{enumerate}
\end{proposition}

\noindent{\it Proof:}
It is easy to check that if $G$ is a modular form then
$G(\tau)-G(\sigma)$ obeys the three-term relations. Conversely,
let us extend $a_{n,m}$ to all integers $n,m$ by setting
$a_{n,m}=0$ if $n$ or $m$ is negative.
The first three-term relation implies that $a_{n-m,m}+a_{n,m-n}
=a_{n,m}$. It follows that $a_{n,m}=a_{n-m,m}$ if $n>m$ and
$a_{n,m}=a_{n,m-n}$ if $m>n$. By the Euclidean algorithm
we see that for $n,m>0$, $a_{n,m}=a_{N,N}$, where $N=(n,m)$
is the greatest common divisor. But $a_{N,N}=0$ since
$Z$ is odd under interchange of $q$ and $r$. Thus only
$a_{n,0}=a_{0,m}$ may be non-zero and 
$Z(\tau,\sigma)=g(\tau)-g(\sigma)$ with $g(\tau)=\sum_{n=1}^\infty a_{n,0}q^n$
In particular, $g(\tau+1)=g(\tau)$. 
The second three-term relation can then be written as
$h(\tau)=h(\sigma)-h(\sigma/\tau)\tau^{-k}$, where $h(\tau)=g(\tau)-\tau^{-k}g(-1/\tau)$
is holomorphic on the upper half plane. If we take the second partial derivative
of this identity with respect
to $\tau$ and $\sigma$, we obtain $z h''(z)+(k+1)h'(z)=0$, with $z=\sigma/\tau$. Thus
$h'$ is homogeneous of degree $-k-1$ and $h(z)=a\tau^{-k}+b$ for some $a$, $b$. Inserting
back in the identity for $h$ shows that $b=-a$. Then $G(\tau)=g(\tau)+a$ obeys
$G(\tau+1)=G(\tau)$ and $G(-1/\tau)=\tau^{k}G(\tau)$ and is thus a modular form of weight $k$.
\hfill $\square$

\section{Elliptic zeta values}
We now define two-parameter deformations of the
values of the zeta function at positive integers,
which we call {\em elliptic zeta values}:
\begin{equation}\label{e-EZV}
Z_k(\tau,\sigma)=
-\frac{(2\pi i)^k}{(k-1)!}
\sum_{j=1}^\infty
j^{k-1}
\frac{q^j-(-1)^kr^j}{(1-q^j)(1-r^j)}\,,\quad
q=e^{2\pi i\tau},\quad r=e^{2\pi i\sigma}, \quad k\in\Z_{\geq1}.
\end{equation}
Clearly, $Z_k(\tau,\sigma)=-(-1)^kZ_k(\sigma,\tau)$.
The relation to the functions $D_k$ is obtained
by expanding the elliptic zeta values 
in a power series in $q$ and $r$. The result is the
following.

If $k$ is even, then
\[
Z_{k}(\tau,\sigma)=D_{k}(r)-D_{k}(q).
\]

 If $k$ is odd, then 
\[
Z_{k}(\tau,\sigma)=
D_{k}(q)+D_{k}(r)+2\sum_{(a,b)=1}D_{k}(q^ar^b).
\]
The sum is over pairs of relatively prime pairs of
positive  integers $a,b$.

\begin{thm}\label{t-1}\ 
\begin{enumerate}
\item[(i)] Let $k\geq4$. Then
$Z_k$ obeys the three-term relations
\begin{eqnarray*}
Z_k(\tau,\sigma)&=&Z_k(\tau,\tau+\sigma)+Z_k(\tau+\sigma,\sigma),\\
Z_k(\tau,\sigma)&=&
\tau^{-k}Z_k\left(-\frac1\tau,\frac\sigma\tau\right)
+(-\sigma)^{-k}Z_k\left(-\frac\tau\sigma,-\frac1\sigma\right),
\end{eqnarray*}
for all $\sigma,\tau\in H$ such that $\sigma/\tau\in H$.
\item[(ii)] $\lim_{\sigma\to 0}\lim_{\tau\to i\infty}
\sigma^kZ_k(\tau,\sigma) =\zeta(k)$ if $k\geq2$.
\item[(iii)]
For $k=1$,  $\lim_{\sigma\to 0}\lim_{\tau\to i\infty}
(\sigma Z_1(\tau,\sigma)+\ln(-2\pi i\sigma))$ is
the Euler constant $\gamma=0.577\dots$ (the logarithm is
real if $\sigma$ is imaginary).
\end{enumerate}
\end{thm}

In particular, by taking the limit of $Z_k$ as $\tau\to i\infty$ we 
obtain
\begin{eqnarray*}
\lim_{\sigma\to0}(\sigma D_1(e^{2\pi i\sigma})+\ln(-2\pi i\sigma))
&=&
\gamma, \\
\lim_{\sigma\to0}\sigma^{k}D_k(e^{2\pi i\sigma})&=&
\zeta(k), \qquad k=2,3\dots
\end{eqnarray*}
We prove Theorem \ref{t-1} in the next section.
\medskip

\noindent{\bf Remark on $\mathrm{SL}(3,\mathbb Z)$.}
The three-term relations in Theorem \ref{t-1}
have an interpretation as $\mathrm{SL}(3,\Z)$
1-cocycle properties parallel
to the ones discussed in \cite{FV} for the
elliptic gamma function. 
We only sketch the construction 
here as the details are the same as in
the case of the elliptic gamma function 
\cite{FV}.
First of all, one notices that the three-term
relations, the symmetry $\tilde Z_k(\tau,\sigma)=
-(-1)^kZ_k(\sigma,\tau)$ and the periodicity
$Z_k(\tau+1,\sigma)=Z_k(\tau,\sigma)$ 
continue to hold if we extend
the domain of $Z_k$ to $\mathrm{Im}\,\tau\neq 0$,
$\mathrm{Im}\,\sigma\neq 0$, by setting
$Z_k(-\tau,\sigma)=Z_k(\tau,-\sigma)=
(-1)^kZ_k(\tau,\sigma)$. We pass to
homogeneous coordinates and set 
$\tilde Z_k(x_1,x_2,x_3)=x_3^{-k}Z_k(x_1/x_3,
x_2/x_3)$. This function is homogeneous
of degree $-k$ and holomorphic on the
dense open set of $\mathbb C^3$ defined
by the conditions $\mathrm{Im}\,x_i/x_3\neq 0$
($i=1,2$). Our identities relate values of
$\tilde Z_k$ at points related by the
 $\mathrm{SL}(3,\mathbb Z)$-action on $\mathbb C^3$: if $k\geq 4$,
\begin{eqnarray*}
\tilde Z_k(x_1,x_2,x_3)&=&\tilde Z_k(x_1+x_3,x_2,x_3)=
\tilde Z_k(x_1,x_2+x_3,x_3)=-\tilde Z_k(x_2,x_1,-x_3),\\
\tilde Z_k(x_1,x_2,x_3)&=&\tilde Z_k(x_1,x_1+x_2,x_3)
+\tilde Z_k(x_1+x_2,x_2,x_3),\\
\tilde Z_k(x_1,x_2,x_3)&=&\tilde Z_k(-x_3,x_2,x_1)+\tilde Z_k(x_1,x_3,-x_2).
\end{eqnarray*}
Let $e_{ij}$ ($i\neq j$) be the elementary
matrices in $\mathrm{SL}(3,\mathbb Z)$, with
$1$ in the diagonal and at the position $(i,j)$,
and zero elsewhere.
These matrices generate $\mathrm{SL}(3,\Z)$.
A consequence of the identities is that
if we set
\[ 
\phi_{e_{12}}(x)=\tilde Z_k(x_1-x_2,x_1,x_3),\qquad
\phi_{e_{32}}(x)=\tilde Z_k(x_2-x_3,x_3,x_1),\qquad
\phi_{e_{ij}}(x)=1,\  j\neq2,
\]
($k\geq4$) then $\phi$ extends (uniquely) to a 1-cocycle
$(\phi_g)_{g\in \mathrm{SL}(3,\mathbb Z)}$
of $\mathrm{SL}(3,\Z)$ with values the
space of holomorphic function with open dense
domain in $\mathbb C^3$, i.e., one has
$\phi_{gh}(x)=\phi_g(x)+\phi_h(g^{-1}x)$ for
all $g,h\in\mathrm{SL}(3,\mathbb Z)$ and $x$
in a dense open set.

\section{The elliptic gamma function}
Theorem \ref{t-1} follows from the results of \cite{FV}
where the properties of modular type for Ruijsenaars's elliptic gamma
function \cite{R}
were discovered. In the normalization of \cite{FV}, the
elliptic gamma function is defined by the double
infinite product
\[
\Gamma(z,\tau,\sigma)=\prod_{j,\ell=0}^\infty
\frac
{1-q^{j+1}r^{\ell+1}e^{-2\pi iz}}
{1-q^{j}  r^{\ell}  e^{ 2\pi iz}}\,\qquad z\in\C,\tau,\sigma\in H.
\]
It obeys the functional relation
$\Gamma(z+\sigma,\tau,\sigma)
=\theta_0(z,\tau)\Gamma(z,\tau,\sigma),
$
which is an elliptic version of Euler's functional equation
for the gamma function. Here $\theta_0$ denotes the modified theta
function
\[
\theta_0(z,\tau)=
\prod_{j=0}^\infty
{(1-q^{j+1}e^{-2\pi iz})}
{(1-q^{j}e^{ 2\pi iz})}\,\qquad z\in\C,\tau\in H.
\]
The Euler gamma function $\Gamma(z)$ is recovered in the limit
\begin{equation}\label{e-SCL}
\Gamma(z)=\lim_{\sigma\to 0}
\lim_{\tau\to i\infty}\theta_0(\sigma,\tau)^{1-z}
\frac
{\Gamma(\sigma z,\tau,\sigma)}
{\Gamma(\sigma  ,\tau,\sigma)}\,.
\end{equation}
\begin{thm}\label{t-2} $($\cite{FV}, p.~54$)$ Suppose 
$\tau,\sigma,\sigma/\tau\in H$. Then
\begin{eqnarray*}
{\Gamma(z,\tau,\sigma)}
&=&
{\Gamma(z+\tau,\tau,\sigma+\tau)}\Gamma(z,\tau+\sigma,\sigma),
\\
\Gamma(
 z/\tau,-1/\tau,\sigma/\tau)&=&
e^{i\pi Q(z;\tau,\sigma)}\label{e-u}
{\Gamma(
{({z-\tau})/\sigma,-\tau/\sigma,-1/\sigma})}
{\Gamma(z,\tau,\sigma)},
\end{eqnarray*}
for some polynomial $Q(z;\tau,\sigma)$ of degree three
in $z$ with coefficients in $\Q(\tau,\sigma)$.
\end{thm}

It follows from the Weierstrass product representation
$\Gamma(z+1)=e^{-\gamma z}
\prod_{j=1}^\infty(1+z/j)^{-1}e^{z/j}$
that the values of $\zeta$ at positive integers are
essentially Taylor coefficients at $1$ of the logarithm
of the Euler gamma function:
\begin{equation}\label{e-logG}
\ln\,\Gamma(z+1)=-\gamma z+\sum_{j=2}^\infty\frac{\zeta(j)}j
(-z)^j.
\end{equation}
Here $\gamma=\lim_{n\to\infty}(\sum_{j=1}^{n-1}1/j-\ln\,n)$ is
the Euler constant.
The elliptic analog of this formula involves the elliptic
zeta values:
\begin{equation}\label{e-logEG}
\ln
\frac
{\Gamma(z+\sigma,\tau,\sigma)}
{\Gamma(\sigma,\tau,\sigma)}=
\sum_{j=1}^\infty \frac{Z_j(\tau,\sigma)}{j}\, (-z)^j.
\end{equation}
This formula, with $Z_j$ defined by eq.~\eqref{e-EZV},
is an easy consequence of the summation formula for
$\ln\,\Gamma(z,\tau,\sigma)$ (\cite{FV} p.~51):
\[
\ln\,\Gamma(z,\tau,\sigma)=
-\frac i2\sum_{j=1}^\infty
\frac{\sin(\pi j(2z-\tau-\sigma))}
{j\sin(\pi j\tau)\sin(\pi j\sigma )}.
\]
{\it Proof of Theorem \ref{t-1}:}\/
The first claim of the theorem 
is proved by taking the
logarithm of the identities of Theorem \ref{t-2} and
expanding them in powers of $z-\sigma$.

 Taking the limit of \eqref{e-logEG}, by using
\eqref{e-SCL} to compare it with \eqref{e-logG} implies
(ii).

In the same way we can deduce (iii).
However, if $n=1$, we have to take into
account the factor $\theta_0^{1-z}$ in
\eqref{e-SCL}: we obtain $\gamma=\lim_{\sigma\to 0}
\lim_{\tau\to\infty}(\sigma Z_1(\tau,\sigma)+
\ln\theta_0(\sigma,\tau))$. 
Since $\ln\theta_0(\sigma,\tau)=\ln(-2\pi i\sigma)$
plus terms that vanish in the limit, the proof is complete.
\hfill$\square$
\medskip

By using the explicit formula for the polynomial $Q$ (see \cite{FV}),
we also obtain the exceptional  relations for $k=1,2,3$:
\begin{thm}\label{t-3}
Let $1\leq k\leq3$ and suppose $\tau,\sigma,\sigma/\tau\in H$. 
Then
\begin{eqnarray*}
Z_k(\tau,\sigma)&=&Z_k(\tau,\tau+\sigma)+Z_k(\tau+\sigma,\sigma),\\
Z_k(\tau,\sigma)&=&
\tau^{-k}Z_k\left(-\frac1\tau,\frac\sigma\tau\right)
+(-\sigma)^{-k}Z_k\left(-\frac\tau\sigma,-\frac1\sigma\right)
+i\pi a_k,
\end{eqnarray*}
where
\[
a_1=
-\frac12
+\frac1{2\tau}
-\frac1{2\sigma}
+\frac{\sigma}{6\tau}
+\frac{\tau}{6\sigma}
+\frac1{6\tau\sigma}\,,
\qquad 
a_2=
-\frac1\tau+\frac1\sigma-\frac1{\tau\sigma}\,,
\qquad 
a_3=
\frac1{\tau\sigma}\,.
\]
\end{thm}

\section{A direct proof of the three-term relations of Theorem \ref{t-1}}\label{s-5}
We thank Don Zagier for providing the following alternative
direct proof of the three term relations of
Theorem \ref{t-1}, which does not use the elliptic gamma
function. As the case where $k$ is even follows easily from
the modular properties of the Eisenstein series, it is sufficient
to consider the case where $k$ is odd.
The proof is based on an alternative formula for $Z_k(\tau,\sigma)$
for $k\geq 5$ odd, which could be seen as a more natural definition of the
elliptic odd zeta values. Let $\epsilon(n)$ be the sign of $n$ 
($\epsilon(n)=1,0,-1$ if $n>0$, $n=0$, $n<0$ respectively) 
and set $\epsilon(a,b)=\frac12(\epsilon(a)+\epsilon(b))$. Then  
\begin{equation}\label{e-alt}
Z_k(\tau,\sigma)=
\sum_{a,b,c}{}'\frac{\epsilon(a,b)}{(a\tau+b\sigma+c)^{k}}, 
\qquad \text{($k$ odd $\geq5$),}
\end{equation}
where $\sum{}'$ means a sum over all $(a,b,c)\neq(0,0,0)$.
To deduce this formula, rewrite \eqref{e-EZV} as
\[
Z_k(\tau,\sigma)=
\frac{(-2\pi i)^k}{(k-1)!}
\sum_{\ell,m=0}^\infty
\sum_{j=1}^\infty
j^{k-1}
({q^j+r^j})q^{j\ell}r^{jm}\,,\quad
q=e^{2\pi i\tau},\quad r=e^{2\pi i\sigma}, \quad \text{$k$ odd,}
\]
then use the Lipschitz formula $\sum_{n\in\mathbb Z}(\rho+n)^{-k}=\frac{(-2\pi i)^k}{(k-1)!}\sum_{j=1}^\infty
j^{k-1}e^{2\pi i\rho j}$, ($\rho\in H$, $k\geq2$).

By using \eqref{e-alt} it is easy to prove Theorem \ref{t-1} (i): the
three term relations follow from the obvious identities
$\epsilon(a-b,b)+\epsilon(a,b-a)=\epsilon(a,b)$,
$\epsilon(a,b)=\epsilon(-c,b)+\epsilon(a,c)$.

It should be possible to prove also the identities in the exceptional
cases where $k=1,3$ using this method, but suitable summation procedures
should be applied to make sense of the series \eqref{e-alt}, which are
not absolutely convergent.

\end{document}